\documentclass{article}
\usepackage{graphicx, amssymb, amsmath}
\usepackage{amsthm}
\usepackage{indentfirst}
\usepackage{mathrsfs}
\usepackage{amsfonts}
\usepackage{color}
\usepackage{graphicx}
\usepackage{algorithm}
\usepackage{algorithmic}
\usepackage{upgreek}
\usepackage{subfigure}

\def\k{\kern .5em}
\def\er{\kern .2em}

\begin{document}
\author{}
\newcommand{\be}{\begin{equation}}
\newcommand{\ee}{\end{equation}}
\newcommand{\ba}{\begin{array}}
\newcommand{\ea}{\end{array}}

\newcommand{\beas}{\begin{eqnarray*}}
\newcommand{\eeas}{\end{eqnarray*}}
\newcommand{\bea}{\begin{eqnarray}}
\newcommand{\eea}{\end{eqnarray}}
\newcommand{\ome}{\Omega}

\newtheorem{theorem}{Theorem}[section]
\newtheorem{lemma}{Lemma}[section]
\newtheorem{remark}{Remark}[section]
\newtheorem{proposition}{Proposition}[section]
\newtheorem{definition}{Definition}[section]
\newtheorem{corollary}{Corollary}[section]

\newtheorem{theo}{Theorem}[section]
\newtheorem{lemm}{Lemma}[section]
\newcommand{\blem}{\begin{lemma}}
\newcommand{\elem}{\end{lemma}}
\newcommand{\bthe}{\begin{theorem}}
\newcommand{\ethe}{\end{theorem}}
\newtheorem{prop}{Proposition}[section]
\newcommand{\bprop}{\begin{proposition}}
\newcommand{\eprop}{\end{proposition}}
\newtheorem{defi}{Definition}[section]
\newtheorem{coro}{Corollary}[section]
\newtheorem{algo}{Algorithm}[section]
\newtheorem{rema}{Remark}[section]
\newtheorem{property}{Property}[section]
\newtheorem{assu}{Assumption}[section]
\newtheorem{exam}{Example}[section]

\renewcommand{\theequation}{\arabic{section}.\arabic{equation}}
\renewcommand{\thetheorem}{\arabic{section}.\arabic{theorem}}
\renewcommand{\thelemma}{\arabic{section}.\arabic{lemma}}
\renewcommand{\theproposition}{\arabic{section}.\arabic{proposition}}
\renewcommand{\thedefinition}{\arabic{section}.\arabic{definition}}
\renewcommand{\thecorollary}{\arabic{section}.\arabic{corollary}}
\renewcommand{\thealgorithm}{\arabic{section}.\arabic{algorithm}}
\newcommand{\lan}{\langle}
\newcommand{\curl}{{\bf curl \;}}
\newcommand{\rot}{{\rm curl}}
\newcommand{\grad}{{\bf grad \;}}
\newcommand{\dvg}{{\rm div \,}}
\newcommand{\ran}{\rangle}
\newcommand{\bR}{\mbox{\bf R}}
\newcommand{\bRn}{{\bf R}^3}
\newcommand{\Coinf}{C_0^{\infty}}
\newcommand{\disp}{\displaystyle}
\newcommand{\ra}{\rightarrow}
\newcommand{\Ra}{\Rightarrow}
\newcommand{\ud}{u_{\delta}}
\newcommand{\Ed}{E_{\delta}}
\newcommand{\Hd}{H_{\delta}}
\newcommand\varep{\varepsilon}
\title{A Decoupling Two-Grid Method for the Steady-State Poisson-Nernst-Planck Equations
}
\author{Ying Yang
\thanks{School of  Mathematics and Computing Science, Guangxi Colleges and Universities Key Laboratory of Data Analysis and Computation, Guangxi Key Laboratory of Cryptography and information Security, Guilin University of Electronic Technology, Guilin, Guangxi 541004, China.
E-mail: yangying@lsec.cc.ac.cn}\and Benzhuo Lu
\thanks{ LSEC,
Institute of Computational Mathematics and Scientific/Engineering
Computing, the National Center for Mathematics and Interdisciplinary
Sciences, Academy of Mathematics and Systems Science,  Chinese
Academy of Sciences, Beijing 100190, China.
E-mail: bzlu@lsec.cc.ac.cn}\and Yan Xie
\thanks{LSEC, Institute of Computational Mathematics and Scientific/Engineering
Computing, Academy of Mathematics and Systems Science, Chinese
Academy of Sciences, Beijing 100190, China}}

\date{}
\maketitle
%\pagenumbering{arabic}

\begin{abstract}
 Poisson-Nernst-Planck equations are widely used to describe the electrodiffusion of ions in a solvated biomolecular system. Two kinds of two-grid finite element algorithms are proposed to decouple the steady-state Poisson-Nernst-Planck equations by coarse grid finite element approximations. Both theoretical analysis and numerical experiments show the efficiency and effectiveness of the two-grid algorithms for solving Poisson-Nernst-Planck equations.
\end{abstract}


\section{Introduction}{\label{sec1}
%part 1. The background of PNP.
%part 2. Some problems in the computing of PNP
%Part 3. Some ways for decoupling the coupled equation
%part 4. The idea of two-grid method and decoupling
Electrodiffusion plays an important role in many fields such as biological ion channels, cellular electrophysiology and semiconductors. For the biological processes, the kinetic properties of them are mainly governed by the electrodiffusion of charged molecules in aqueous solution. The numerical methods for deriving the kinetic parameters usually include discrete methods (such as Monte Carlo, Brownian
dynamics and Langevin dynamics) and continuum methods. The latter is more efficient for simulating large systems and easier to be modified to include more physical functions. The electrodiffusion processes in biomolecular systems are usually described by a continuum model called Poisson-Nernst-Planck (PNP) equations, which is regarded as one of the most efficient theoretical methods for studying electrodiffusion.

The PNP equations  are a coupled system of nonlinear partial differential equations consisting of the Nernst-Planck equation and the electrostatic Poisson equation. The steady-state PNP equations in the biomolecular modeling are as follows \cite{luho10}:
 \bea
\label{eqn0}\left\{\begin{array}{ll}\nabla\cdot D^i\left(\nabla
p^i+\beta q^ip^i\nabla\phi\right)=0,~~\mbox{in}~~ \ome_s,~1\leq i\leq n,\\
-\nabla\cdot(\epsilon\nabla\phi)-\lambda\displaystyle{\sum_{i=1}^n}
q^ip^i=\rho^f,~~\mbox{in}~~\ome,\end{array}\right.\eea where $\phi$
is the electrostatic potential and $p^i$ is the concentration of the
$i$-th ion species. This model is used to describe the
electrodiffusion of mobile ions in a solvated biomolecular system (see Fig. 1). %Fig. 1 illustrates a solvated biomolecular system in an
%open domain $\ome\subset \mathbb{R}^3$.
Since the analytic solutions of the PNP equations only exit in very few cases for simple shape molecules, a variety of numerical methods have been proposed to solve them \cite {caco00,eich93,luho10,wusr02,zhlu08}. For example, the finite difference method has been widely used to solve the PNP equations describing electrodiffusion in biological ion channels or other transmembrane pores \cite{bosa09,imro02,kuco99,zhch11}, but the accuracy is not so high when it is applied to the biomolecular models with highly irregular surfaces. The finite element method is considered to be very promising in which irregular shapes can be fitted more easily when applying to the PNP equations. \cite {jeke91,luho10,luzh07,sozh041,sozh042,yalu13}.

\begin{figure}[htbp]
  \subfigure[]{
    %\label{fig:subfig:b} %% label for second subfigure
    %\includegraphics[width=4.5cm,totalheight=4.5cm]{illu.eps}}
    \includegraphics[width=4.5cm,totalheight=4.5cm]{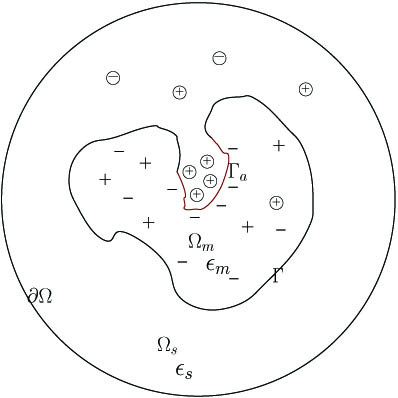}}
    \centering
     %\subfigure[]{
    %\label{fig:subfig:c} %% label for second subfigure
   % \includegraphics[width = 3.5cm, totalheight = 4.5cm]{diALA_S.eps}}
 \caption{2-D illustration of the computational
domain modeling a solvated biomolecular system. The domain $\ome_m\subset\ome$ represents the
biomolecule(s) and the remain domain $\ome_s=\ome\setminus \bar{\ome}_m$
shows a solvent surrounding the biomolecule(s). The molecular
surface $\Gamma$ interfaces domains $\ome_m$ and $\ome_s$. Charged
ligands in this model are also treated as diffusive species, and
might react with the biomolecules on a part of the molecular surface
$\Gamma_a$. The diffusive particles are distributed in
$\ome_s$. }
\label{fig:biomolecular}
\end{figure}

In general, there are two types of approaches to solve such a multimodel problem like PNP equations. One is to consider the equations as a large system and solve the overall system together. The other is to first decouple the system and then solve the equations respectively. The advantages of the latter are obvious. For example, it can be implemented more easily and efficiently since the computational systems are smaller; it can effectively exploit the existing computing softwares and it can result in parallelism in some cases. The Gummel iteration \cite{gu64} is such a type of decoupling approach which may be the most commonly used decoupling method for solving PNP equations \cite{bu11,je87,luzh11}.
This approach is to solve one of the equations first and then substitute the solution in successive equations. For example, consider the following system coupled by two equations
$$ F_1(u_1,u_2)=0, ~~F_2(u_1,u_2)=0.$$
The Gummel iteration for the above system could be: given $u_2^0$, for $k\ge 0$, find ($u_1^{k+1},u_2^{k+1}$) such that
$$ F_1(u_1^{k+1},u_2^k)=0,~~F_2(u_1^{k+1},u_2^{k+1})=0,$$
until the error between the $(k+1)$th solution and $k$th solution is less than the tolerance.
%The Gummel iteration \cite{gu64} can be considered as a special type of block G-S iteration and is wide used to solving PNP equations for semiconductor device.
However, the Gummel iteration converges slowly or even diverges for PNP equations modeling complex biomoleculars such as protein and DNA. There are some other approaches to decouple multimodel problems. Most of them are based on the idea of domain decomposition \cite{camu09,glpa97,maho00,quva99}. Among those approaches, the two-grid finite element method closely combines with the finite element method. The two-grid finite element method was originally proposed by Xu for the partial differential equation (PDE) to deal with the asymmetry and indefiniteness \cite{xu92,xu94,xu96} and gradually developed for some other applications such as linearization, localization and parallelization \cite{maxu95,xuzh00,xuzh02}.
%All of these designs of two-grid finite element method are mainly used for a single PDE.
Compared with the work for the single PDE as above, the primary motivation of the two-grid finite element method for the coupled PDEs (which involves two or more equations) is different. It can be applied to decouple the system of PDEs and has shown the efficiency and effectiveness for some coupled systems \cite{camu09,jish06,muxu07}. Jin and his coauthors \cite {jish06} used the
two-grid discretization method to decouple the Schr$\ddot{o}$dinger equation arising from quantum mechanics. Mu and Xu \cite{muxu07} presented a type of decoupling method based on two-grid finite element method and applied it to solving Stoke-Darcy model for coupling fluid flow with porous media flow.

   In this paper, we will use the two-grid finite element method to decouple the system of steady-state PNP equations. Since the equations are quite different from the above multimodel problems, the design of the two-grid method could not be directly inspired from the above work and the resulting two-grid schemes are also different. We shall design two decoupling schemes by using the two-grid method. One is semi-decoupled. The other is fully decoupled and is suitable for parallelism. These schemes can provide good initial values for solving PNP equations and do not require iterations on the fine discrete grid as the Gummel iteration does, which can improve the computational efficiency and save the computational time. We also get some
analysis results. The results show that if the finite element solution on the
coarse grid approximates that on the fine grid well enough, then the two-grid method can achieve the similar approximation effect that the conventional finite element method could do. %
Some numerical examples including a biomolecular problem are shown to verify the theoretic results.%

The paper is organized as follows. In the following section, some preliminaries including the introduction to the PNP equations are presented. In Section 3, two decoupling two-grid algorithms are proposed and the corresponding error analyses are presented. Some numerical examples are shown in Section 4.

%---------------------------

\setcounter{equation}{0}
\section{Preliminaries}\label{sec2}

 In this section, we shall first describe the steady-state
PNP equations and boundary conditions, then introduce the corresponding weak formulations. The finite element approximations are also studied in
this section.

\par In this paper, we consider the following steady-state PNP system for simulating biomolecular diffusion-reaction process \cite{luho10,yalu13}

\bea\label{prob1} \left\{\begin{array}{ll}\nabla\cdot
D^i\left(\nabla
p^i+\beta q^ip^i\nabla\phi\right)=0,~~\mbox{in}~~ \ome_s,~1\leq i\leq n,\\
-\nabla\cdot(\epsilon\nabla\phi)-\lambda\displaystyle{\sum_{i=1}^n}
q^ip^i=\rho^f,~~\mbox{in}~~\ome\subset R^3\end{array}\right.\eea %%
with the following interface conditions and boundary conditions (for
simplicity, the reactive molecular surface $\Gamma_a$ is not
considered, see Fig. 1) \bea\label{bd1}\left\{\begin{array}{ll}
[\phi]_{\Gamma}=0,~~[\epsilon\frac{\partial\phi}{\partial
n}]_{\Gamma}=0,~~\mbox{on}~~\Gamma=\bar{\ome}_s\bigcap\bar{\ome}_m,\\
n\cdot D^i\left(\nabla p^i+\beta q^i
p^i\nabla\phi\right)=0,~~\mbox{on}~~\Gamma,~1\leq i\leq n,\\
\phi=0,~~\mbox{on}~~\partial\ome,\\
p^i=p_{pulk}^i,~~\mbox{on}~~\partial\ome, \end{array}\right.\eea
 where $p^i(x)$
is the concentration of the $i$-th species particle carrying charge
$q^i$, $1\leq i\leq n$, $\phi(x)$ is the electrostatic potential, $D^i(x)$ is the
diffusion coefficient, the permanent (fixed) charge distribution $\rho^f(x)=\displaystyle\sum_{j=1}^{N_m}
q_j\delta(x-x_j)$ is a linear combination of Dirac Delta functions and represents an ensemble of singular charges $q_j$ located
at $x_j$ inside biomolecules, $\beta=1/(\kappa_BT)$ is the inverse
Boltzmann energy, $\epsilon(x)=\left\{\begin{array}{ll}\epsilon_m,
~~ x\in\ome_m,\\
\epsilon_s,~~x\in\ome_s,\end{array}\right.$ is the dielectric
coefficient, $\lambda=\left\{\begin{array}{ll}0,
~~ \mbox{in}~~\ome_m,\\
1,~~\mbox{in}~~\ome_s\end{array}\right.$, $n$ is the outer normal vector and
$p_{pulk}^i,~i=1,2\cdots,n$ are given functions.

 Let
$\ome\subset \mathbb{R}^3$ be a polyhedral convex domain with a
Lipschitz-continuous boundary $\partial\ome$. Assume the
interface $\Gamma$ is sufficiently smooth, say, of class $C^2$. We shall adopt the
standard notations for Sobolev spaces $W^{s,p}(\Omega)$ and their
associated norms and seminorms, see, e.g., \cite{ad75, brsc}. For
$p=2$, we denote $H^s(\Omega)=W^{s,2}(\Omega)$ and
$H^1_0(\Omega)=\{v|v\in H^1(\Omega): v\mid_{\partial\Omega}=0\}$,
where $v\mid_{\partial\Omega}=0$ is in the sense of trace,
$\|\cdot\|_{s,p,\Omega}= \|\cdot\|_{W^{s,p}(\Omega)}$ and
$(\cdot,\cdot)$
is the standard $L^2$-inner product. The weak formulations of (\ref{prob1}) and (\ref{bd1}) are that:
Find $\phi\in H^1(\ome_s)$
and $p^i\in V=\{v|v\in H^1(\ome_s),v|_{\partial\ome}=p^i_{pulk}\}$
$(1\leq i\leq n)$ such that \cite{yalu13}%%
 \bea\label{wf1}a_1(p^i,v)+b_1(p^i,\phi,v)=0,~~\forall v\in
V_0,\eea \bea\label{wf2}a_2(\phi,w)+b_2(p^i,w)=f(w),~~\forall
w\in C_0^{\infty}(\ome),\eea where \bea
\label{a1b1}a_1(p^i,v)=(D^i\nabla p^i,\nabla
v),~~~~b_1(p^i,\phi,v)=(D^i\beta q^i p^i\nabla \phi,\nabla v), \eea
\bea \label{a2b2} a_2(\phi,w)=(\epsilon \nabla \phi,\nabla
w),~~~~b_2(p^i,w)=(-\lambda\sum_{i=1}^n  q^i p^i,w),~~f(w)=(\rho^f,w), \eea and the space
$V_0=\{v|v\in H^1(\ome_s),v|_{\partial\ome}=0\}$. Here $D^i(x)\geq D_0>0$, $\epsilon(x),~\beta,~D^i(x)$ and $q^i\in L^{\infty}(\ome),~(1\leq i\leq n)$.

Assume that $T^h(\ome)$ is a quasi-uniform mesh of size $h\ll 1$ . For ease of analysis, we suppose that the triangulation resolves the interface,
although this assumption may be weakened in the practical computation.
We define linear finite element spaces
\bea\label{femspace} S^h(\ome)=\{v\in H^1(\ome): v|_{e}\in P^1(e), \forall e
\in T^h(\ome)\},
   ~~~~S_0^h(\ome)=S^h(\ome) \cap H_0^1(\Omega),\eea
where $P^1(e)$ is the set of linear polynomials. The coarse spaces $S^H(\ome)$ and $S_0^H(\ome)$ are defined by replacing $h$ with $H$ in (\ref{femspace}).
\par Suppose there exists a unique solution ($\phi$,~$p^i$) satisfying (\ref{wf1})-(\ref{a2b2}). The standard
finite element discrete scheme for (\ref{wf1})-(\ref{a2b2})
reads:
\begin{algo}\label{algo0}(Standard finite element method \cite{yalu13})
 Find $\phi_h\in S_0^h(\ome)$
and $p_h^i\in S^h(\ome_s)\cap V,~(1\leq i\leq n)$, such that
 \bea\label{fea1}a_1(p^i_h,v_h)+b_1(p^i_h,\phi_h,v_h)=0,~~\forall v_h\in
V_0\cap S^h(\ome_s),~~1\leq i\leq n,\eea
\bea\label{fea2}a_2(\phi_h,w_h)+b_2(p^i_h,w_h)=\hat{f}(w_h),~~\forall
w_h\in S_0^h(\ome),\eea
where $\hat{f}(w_h)=\displaystyle{\sum_{j=1}^{N_m}}q_jw_h(x_j)$ is an approximation to the functional $f(w)$ \cite{yazh06}.
\end{algo}
We assume there exists a unique solution
($\phi_h$,~$p^i_h$) satisfying (\ref{fea1}) and (\ref{fea2}). Some error bounds were presented in \cite{yalu13} for the finite element approximation. For example, if $\phi\in H^{1+m}(\ome_s) $ and $p^i\in
H^{1+m}(\ome_s)~(1\leq i\leq n),~0<m\leq 1$ , then we have
\bea
\label{phi1os}\|\phi-\phi_h\|_{1,\ome_s}\leq C
(h^m+\sum_{i=1}^n\|p^i-p_h^i\|_{0,\ome_s}),\eea
and \bea \label{errp}
\|p^i-p_h^i\|_{1,\ome_s}\leq C(h^m+\sum_{i=1}^n
\|p^i-p^i_h\|_{0,\ome_s}),\eea when $p_h^i\in L^{\infty}(\ome_s)$.
Although there is no error estimate in $L^2$ norm for the solution $p_h^i$ of the steady-state PNP equations, the numerical results show that the second order accuracy could be achieved (see \cite{yalu13} and also the results in Table $2$ of Section $4$ in this paper).

We also note that the wellposedness and error estimations in $H^1$ and $L^2$ norms of the finite element approximation for the time-dependent PNP are presented in \cite{prsc09} and \cite{susu16}, respectively.

We introduce the auxiliary problem which shall be used in the next section: Find $w\in H^1(\ome_s)$, such that
\bea\label{dualprob} a_1(v,w)+b_1(v,\phi,w)=(f,v),~~\forall v\in H_0^1(\ome_s),\eea
where $f\in L^2(\ome_s)$. The finite element discrete scheme for (\ref{dualprob}) reads: Find $w_h\in S^h(\ome_s)$, such that
\bea \label{dualapp}a_1(v_h,w_h)+b_1(v_h,\phi,w_h)=(f,v_h),~~\forall v_h\in S_0^h(\ome_s).\eea
\begin{lemm}\label{reglemm}If there exists a unique solution for the problem (\ref{dualprob}) when $f\in L^2(\ome_s)$, then the following regularity result holds (cf. Grisvard \cite{gr85})
\beas\label{regres} \|w\|_{2,\ome_s}\leq C\|f\|_{0,\ome_s}.\eeas
Furthermore, if $\phi\in W^{2,\infty}(\ome_s)$, then (cf. Xu \cite{xuzh00})
\beas\label{sta}\|\nabla w_h\|_{0,\ome_s}\leq C \|\nabla w\|_{0,\ome_s}.\eeas
\end{lemm}

 %%%%%%%%%%%%%%%%%%%%%%%%%%%%%%%%%%%%%%%%%%%%%%%%%%%%%%%%%%%%%%%%%%%%%%%%%%%%%5
 \section{The two-grid finite element method}\label{sec3}\setcounter{equation}{0}
In this section, we shall present the two-grid finite element method
for the PNP equations. Two algorithms are provided to decouple the
strong coupled equations. The first one is a semi-decoupling scheme. The second one is a fully decoupling scheme, which is suitable for parallel computing. Some error analyses are also derived for the
two-grid finite element approximations.

\begin{algo}\label{algo1}(Two-grid algorithm I)$     $\\
Step 1. Solve the coupled problem on the coarse grid: find $p_H^i\in
S^H(\ome_s)\cap V~(1\leq i\leq n)$ and $\phi_H\in S_0^H(\ome)$, such that \bea
\label{coapi}a_1(p_H^i,v_H)+b_1(p_H^i,\phi_H,v_H)=0,~~\forall v_H\in
S^H(\ome_s)\cap V_0,~~1\leq i\leq n,\eea
\bea\label{coaphi}a_2(\phi_H,w_H)+b_2(p_H^i,w_H)=\hat{f}(w_H),~~\forall w_H\in S_0^H(\ome),\eea
where $\hat{f}(w_H)=\displaystyle{\sum_{j=1}^{N_m}}q_jw_H(x_j)$.\\
Step 2. We first solve the Poisson equation on the fine grid: find $\phi_h^*\in
S_0^h(\ome)$, such that
\bea\label{finphi}a_2(\phi_h^*,w_h)+b_2(p_H^{i},w_h)=\hat{f}(w_h),~~\forall w_h\in S_0^h(\ome),\eea
then we solve the
Nernst-Planck equation on the fine grid: find $p_h^{i,*}\in
S^h(\ome_s)\cap V$, such that \bea
\label{finpi}a_1(p_h^{i,*},v_h)+b_1(p_h^{i,*},\phi_h^*,v_h)=0,~~\forall
v_h\in S^h(\ome_s)\cap V_0,\eea
where $\hat{f}(w_h)=\displaystyle{\sum_{j=1}^{N_m}}q_jw_h(x_j)$.
%\bea\label{finphi}a_2(\phi_h^*,w_h)+b_2(p_H^{i},w_h)=\sum_j
%q_jw_h(x_j),~~\forall w_h\in S_0^h(\ome).\eea
\end{algo}
Obviously, (\ref{coapi}) and (\ref{coaphi}) in Step 1 in fact are the standard finite element discretization on a coarse space. Hence, if the Gummel iteration %(\ref{gseqn1})-(\ref{gseqn2}) (in which $h$ is replaced by $H$)
is applied in this step, the iteration between the coupled equations is still required but it converges faster because of much less of the degree of freedoms comparing with Algorithm \ref{algo0}.
The system in Step 2 is a semi-decoupling one. To solve
this system, the solution $\phi_h^*$ in (\ref{finphi}) should be
solved first and then be inserted into (\ref{finpi}) to get the
solution $p_h^{i,*}$. Both (\ref{finphi}) and (\ref{finpi}) need to be solved only for one time, respectively. Hence this algorithm could reduce plenty of computational time comparing with Algorithm \ref{algo0} by using the Gummel iteration which requires a lot of iterations between the equations solving individually. Moreover, it can naturally avoid the slow convergence or divergence of the Gummel iteration.

The following two theorems provide the error bounds for the two-grid
solutions of Algorithm \ref{algo1} and the solutions of
Algorithm \ref{algo0}.

\begin{theo} \label{theo1}If $(\phi_h,p_h^i)$, $(\phi_H,p_H^i)$ and $(\phi_h^*,p_h^{i,*})$ are the solutions of (\ref{fea1})-(\ref{fea2}), (\ref{coapi})-(\ref{coaphi})
 and (\ref{finphi})-(\ref{finpi}), respectively, then
\bea \label{theo1eqn}\|\nabla(\phi_h-\phi_h^*)\|_{0,\ome}\leq C
\sum_{i=1}^n\|p_h^i-p_H^i\|_{0,\ome_s}.\eea
\end{theo}
\textit{Proof.} From (\ref{fea2}) and (\ref{finphi}), we have \beas
a_2(\phi_h-\phi_h^*,w_h)+b_2(p_h^i-p_H^i,w_h)=0,~~\forall w_h\in
S_0^h(\ome).\eeas Taking $w_h=\phi_h-\phi_h^*$ in the above
equality, we obtain
\beas \|\nabla(\phi_h-\phi_h^*)\|_{0,\ome}^2&\leq& |a_2(\phi_h-\phi_h^*,\phi_h-\phi_h^*)|\\
         &=&|b_2(p_h^i-p_H^i,\phi_h-\phi_h^*)| \\
         &\leq & \sum_{i=1}^n\|p_h^i-p_H^i\|_{0,\ome_s} \|\phi_h-\phi_h^*\|_{0,\ome_s}. \eeas
Hence, \beas \|\nabla(\phi_h-\phi_h^*)\|_{0,\ome}\leq C
\sum_{i=1}^n\|p_h^i-p_H^i\|_{0,\ome_s}.\eeas This completes the
proof.$\hfill\Box$

The following error bounds in $L^2$ norm will be used for presenting the error estimate for the concentration in $H^1$ norm.
\begin{lemm}\label{lemm1} Suppose the assumptions of Lemma \ref{reglemm} and Theorem \ref{theo1} hold, $p^{i,*}_h\in L^{\infty}(\ome_s)$, $\phi$ is the solution of (\ref{wf1})-(\ref{wf2}) and $\phi_h$ approximates $\phi$ well enough, then we have
\beas\|p_h^i-p_h^{i,*}\|_{0,\ome_s}\leq
C\sum_{i=1}^n\|p_h^i-p_H^i\|_{0,\ome_s}.\eeas
\end{lemm}
\textit{Proof.} Set $e=p_h^i-p_h^{i,*}$. Let $w$ be the solution of the following auxiliary problem:

\beas\label{theo2eqn0} a_1(v,w)+b_1(v,\phi,w)=(e,v),~~ \forall v\in
H^1_0(\ome_s),\eeas
and $w_h$ be the finite element approximation to $w$ satisfying
\bea \label{femwh} a_1(v_h,w_h)+b_1(v_h,\phi,w_h)=(e,v_h),~~ \forall v_h\in
S^h_0(\ome_s).\eea

Taking $v_h=e$ in (\ref{femwh}), then from (\ref{fea1}) and (\ref{finpi}), we
have
\bea\label{theo2eqn1}\|e\|_{0,\ome_s}^2
%&=&a_1(e,w_h)+b_1(e,\phi,w_h)\nonumber \\
                              &=&b_1(p_h^{i,*},\phi_h^*,w_h)-b_1(p_h^{i},\phi_h,w_h)+b_1(p_h^{i}-p_h^{i,*},\phi,w_h)
                              \nonumber\\
                              %&=&b_1(p_h^{i,*},\phi_h^*-\phi,w_h)+b_1(p_h^{i},\phi-\phi_h,w_h)
                              %\nonumber\\
                              &=&b_1(p_h^{i,*},\phi_h^*-\phi_h+\phi_h-\phi,w_h)+b_1(p_h^{i},\phi-\phi_h,w_h)\nonumber\\
                              &=&b_1(p_h^{i,*},\phi_h^*-\phi_h,w_h)+b_1(p_h^{i,*}-p_h^{i},\phi_h-\phi,w_h)\nonumber\\
                              &\leq& C(\|\nabla(\phi_h-\phi_h^*)\|_{0,\ome_s}+\|p_h^{i,*}-p_h^i\|_{0,\ome_s}\|\nabla(\phi-\phi_h)\|_{0,\infty,\ome_s}) \|\nabla w_h\|_{0,\ome_s},
%                              &=&(I)_1+(I)_2,
\eea
where the assumption $p^{i,*}_h\in L^{\infty}(\ome_s)$ is used. From Theorem \ref{theo1} and Lemma \ref{reglemm}, we get
\beas \|e\|_{0,\ome_s}^2\leq C(\sum_{i=1}^n\|p_h^i-p_H^i\|_{0,\ome_s}+\|e\|_{0,\ome_s}\|\nabla(\phi-\phi_h)\|_{0,\infty,\ome_s})\|e\|_{0,\ome_s}.\eeas
If $\phi_h$ approximates $\phi$ well enough satisfying $C\|\nabla(\phi-\phi_h)\|_{0,\infty,\ome_s}<1$, then we
can easily get the result of Lemma \ref{lemm1}. $\hfill\Box$

\begin{theo}\label{theo3}Suppose the assumptions of Lemma \ref{lemm1} hold and $\phi_h\in W^{1,\infty}(\ome_s)$,
then we have
\beas\|\nabla(p_h^i-p_h^{i,*})\|_{0,\ome_s}\leq
C\sum_{i=1}^n\|p_h^i-p_H^i\|_{0,\ome_s}.\eeas
\end{theo}
\textit{Proof.} Similarly, set $e=p_h^i-p_h^{i,*}$. From (\ref{fea1}) and (\ref{finpi}), we have
\beas a_1(e,v_h)+b_1(p_h^i,\phi_h,v_h)-b_1(p_h^{i,*},\phi_h^*,v_h)=0,~~\forall v_h\in S^h(\ome_s)\cap V_0.\eeas
\bea \label{theo3eqn1}\|\nabla e\|_{0,\ome_s}^2&\leq& C |a_1(e,e)|\nonumber\\
&=&C|b_1(p_h^i,\phi_h,e)-b_1(p_h^{i,*},\phi_h^*,e)|\nonumber\\
&\leq& C \|p_h^i\nabla\phi_h-p_h^{i,*}\nabla\phi_h^*
\|_{0,\ome_s}\|\nabla e\|_{0,\ome_s}.\eea If $p^{i,*}_h\in
L^{\infty}(\ome_s)$ and $\phi_h\in W^{1,\infty}(\ome_s)$, then
\bea\label{theo3eqn2} \|p_h^i\nabla\phi_h-p_h^{i,*}\nabla\phi_h^*) \|_{0,\ome_s}&=&\|(p_h^i-p_h^{i,*})\nabla\phi_h+p_h^{i,*}(\nabla\phi_h-\nabla\phi_h^*)\|_{0,\ome_s}\nonumber\\
&\leq&
\|e\|_{0,\ome_s}+\|\nabla\phi_h-\nabla\phi_h^*\|_{0,\ome_s}.\eea
Inserting (\ref{theo3eqn2}) into (\ref{theo3eqn1}) and using
Theorem \ref{theo1} and Lemma \ref{lemm1}, we can get the result of Theorem \ref{theo3}.$\hfill\Box$

%(\ref{theo1eqn}), we have \beas\|\nabla e\|_{0,\ome_s}^2&\leq&
%C(\|e\|_{0,\ome_s}+\sum_{i=1}^n\|p_h^i-p_H^i\|_{0,\ome_s})\|\nabla e\|_{0,\ome_s}\\
%&\leq& C (\sum_{i=1}^n\|p_h^i-p_H^i\|_{0,\ome_s}+h\|\nabla e\|_{0,\ome_s})\|\nabla e\|_{0,\ome_s}\\
%&\leq& C_1 \sum_{i=1}^n\|p_h^i-p_H^i\|_{0,\ome_s}^2+C_2\|\nabla e\|_{0,\ome_s}^2+Ch\|\nabla e\|_{0,\ome_s}^2,\eeas
%where we have used Lemma 3.1, $C_1$ and $C_2$ are two constants with $C_2<1$. If the mesh size $h$ is small enough satisfying $Ch<1$, then we
%immediately complete the proof.$\hfill\Box$

\begin{remark} \label{rema1}Theorem \ref{theo1} and Theorem \ref{theo3} show that
the errors between the finite element solution $(\phi_h,p_h^i)$ in Algorithm \ref{algo0} and
the two-grid solution $(\phi_h^*,p_h^{i,*})$ in Algorithm \ref{algo1} are
controlled by the error of the solution for the concentration on the
coarse grid and on the fine grid. If the concentration on the coarse grid
approximates to that on the fine
grid well enough, %then the error
%$\sum_{i=1}^n\|p_h^i-p_H^i\|_{0,\ome_s}$ may achieve optimal convergence
%rate, i.e.
for example,
\beas \sum_{i=1}^n\|p_h^i-p_H^i\|_{0,\ome_s}=O(H^2),\eeas
then from Theorem \ref{theo1} and Theorem
\ref{theo3} we have \bea
\|\nabla(\phi_h-\phi_h^*)\|_{0,\ome}=O(H^2)=O(h)\eea and \bea \|
\nabla(p_h^i-p_h^{i,*})\|_{0,\ome_s} =O(H^2)=O(h),\eea when $h=H^2$. This
means the two-grid method can achieve the similar effect as the
classic finite element method could do (since the optimal convergence rate
for classic finite element method in $H^1$ norm is not better than
O(h), cf. (\ref{phi1os}) and (\ref{errp})), if the solution for the concentration on the coarse grid approximates
that on the fine grid well enough.
\end{remark}

%%%%%%%%%%%%%%%%%%%%%%%%%%%%%%%%%%%%%%%
Next, we shall present a fully decoupling two-grid algorithm.
\begin{algo}\label{algo2}(Two-grid algorithm II)$     $\\
Step 1. Solve the coupled problem on the coarse grid: find $p_H^i\in
S^H(\ome_s)\cap V~(1\leq i\leq n)$ and $\phi_H\in S_0^H(\ome)$, such that \beas
\label{coapi2}a_1(p_H^i,v_H)+b_1(p_H^i,\phi_H,v_H)=0,~~\forall
v_H\in S^H(\ome_s)\cap V_0,~~1\leq i\leq n,\eeas
\beas\label{coaphi2}a_2(\phi_H,w_H)+b_2(p_H^i,w_H)=\hat{f}(w_H),~~\forall w_H\in S_0^H(\ome).\eeas
 Step 2. Solve the Nernst-Planck equation and Poisson equation: find
$p_h^{i,*}\in S^h(\ome_s)\cap V~(1\leq i\leq n)$ and $\phi_h^*\in S_0^h(\ome)$, such
that \bea
\label{finpi2}a_1(p_h^{i,*},v_h)+b_1(p_h^{i,*},\phi_H,v_h)=0,~~\forall
v_h\in S^h(\ome_s)\cap V_0,~~1\leq i\leq n,\eea
\bea\label{finphi2}a_2(\phi_h^*,w_h)+b_2(p_H^{i},w_h)=\hat{f}(w_h),~~\forall w_h\in S_0^h(\ome).\eea
\end{algo}
%Comparing with Algorithm \ref{algo1}, another finite element solution $\phi_H$ besides $p_H^i$ is also used to decouple the system on the
%fine grid.
In this algorithm, two coarse grid finite element approximations $\phi_H$ and $p_H^i$ are used to decouple the system on the fine space, which is different from the case in Algorithm 3.1.
The system (\ref{finpi2})-(\ref{finphi2}) are fully
decoupled, so it can be solved in parallel on the fine grid level.

\par The following are the error analysis for Algorithm
\ref{algo2}. First, for Algorithm \ref{algo2}, we have
\begin{theo} \label{theo4} If $(\phi_h,p_h^i)$, $(\phi_H,p_H^i)$ and $(\phi_h^*,p_h^{i,*})$ are the solutions of (\ref{fea1})-(\ref{fea2}), (\ref{coapi})-(\ref{coaphi})
 and (\ref{finpi2})-(\ref{finphi2}), respectively, then
 %The following estimate holds under the assumption of Theorem \ref{theo1},
% If $\phi_h$ and $\phi_h^*$ are the solutions of (\ref{fea1})-(\ref{fea2}) and (\ref{finpi2})-(\ref{finphi2}), respectively, then
\bea \label{theo4eqn0}\|\nabla(\phi_h-\phi_h^*)\|_{0,\ome}\leq C
\sum_{i=1}^n\|p_h^i-p_H^i\|_{0,\ome_s}.\eea
\end{theo}
The proof is the same as Theorem \ref{theo1}, since the only
difference between Algorithm \ref{algo1} and \ref{algo2} is (\ref{finpi2}) which
is not used in this proof.

\begin{lemm}\label{lemm3} Suppose the assumptions of Lemma \ref{reglemm} and Theorem \ref{theo4} hold, $p^{i,*}_h\in L^{\infty}(\ome_s)$, $\phi$ is the solution of (\ref{wf1})-(\ref{wf2}) and $\phi_h$ approximates $\phi$ well enough, then we have
\beas\|p_h^i-p_h^{i,*}\|_{0,\ome_s}\leq
C\|\nabla\phi_h-\nabla\phi_H\|_{0,\ome_s}.\eeas
\end{lemm}
\textit{Proof.} The proof is similar to that of Lemma \ref{lemm1}.  Let $e=p_h^i-p_h^{i,*}$. From (\ref{fea1}) and (\ref{finpi2}), we have
\beas a_1(e,w_h)=b_1(p_h^{i,*},\phi_H,w_h)-b_1(p_h^{i},\phi_h,w_h).\eeas
Taking $e=p_h^i-p_h^{i,*}$ in (\ref{femwh}) and from the above, we deduce that
\beas\label{lemm2eqn1}\|e\|_0^2&=&b_1(p_h^{i,*},\phi_H,w_h)-b_1(p_h^{i},\phi_h,w_h)+b_1(p_h^{i}-p_h^{i,*},\phi,w_h)\nonumber \\
                              &=&b_1(p_h^{i,*},\phi_H-\phi_h+\phi_h-\phi,w_h)+b_1(p_h^{i},\phi-\phi_h,w_h)\nonumber\\
                              &=&b_1(p_h^{i,*},\phi_H-\phi_h,w_h)+b_1(p_h^{i,*}-p_h^{i},\phi_h-\phi,w_h)\nonumber\\
                              &\leq& C(\|\nabla(\phi_h-\phi_H)\|_{0,\ome_s}+\|p_h^{i,*}-p_h^i\|_{0,\ome_s}\|\nabla(\phi-\phi_h)\|_{0,\infty,\ome_s}) \|\nabla w_h\|_{0,\ome_s},\nonumber \eeas

where the assumption $p^{i,*}_h\in L^{\infty}(\ome_s)$ is used. From Lemma \ref{reglemm} and the assumption $\phi_h$ approximates $\phi$ well enough satisfying $C\|\nabla(\phi-\phi_h)\|_{0,\infty,\ome_s}<1$, we can obtain the result of Lemma \ref{lemm3}. $\hfill\Box$

\begin{theo}\label{theo5} Suppose the assumptions of Lemma \ref{lemm3} hold and $\phi_h\in W^{1,\infty}(\ome_s)$,,
 then we have
\beas\|\nabla(p_h^i-p_h^{i,*})\|_{0,\ome_s}\leq
C\|\nabla \phi_h-\nabla \phi_H\|_{0,\ome_s}.\eeas
\end{theo}
\textit{Proof.} Set $e=p_h^i-p_h^{i,*}$. From (\ref{fea1}) and
(\ref{finpi2}), we have \beas
a_1(e,v_h)+b_1(p_h^i,\phi_h,v_h)-b_1(p_h^{i,*},\phi_H,v_h)=0,~~\forall
v_h\in S^h(\ome_s)\cap V_0.\eeas
\bea \label{theo4eqn1}\|\nabla e\|_{0,\ome_s}^2&\leq& C |a_1(e,e)|\nonumber\\
&=&C|b_1(p_h^i,\phi_h,e)-b_1(p_h^{i,*},\phi_H,e)|\nonumber\\
&\leq& C \|p_h^i\nabla\phi_h-p_h^{i,*}\nabla\phi_H
\|_{0,\ome_s}\|\nabla e\|_{0,\ome_s}.\eea If $p^{i,*}_h\in
L^{\infty}(\ome_s)$ and $\phi_h\in W^{1,\infty}(\ome_s)$, then
\bea\label{theo4eqn2} \|p_h^i\nabla\phi_h-p_h^{i,*}\nabla\phi_H \|_{0,\ome_s}&=&\|(p_h^i-p_h^{i,*})\nabla\phi_h+p_h^{i,*}(\nabla\phi_h-\nabla\phi_H)\|_{0,\ome_s}\nonumber\\
&\leq&
\|e\|_{0,\ome_s}+\|\nabla\phi_h-\nabla\phi_H\|_{0,\ome_s}.\eea
The proof of Theorem \ref{theo5} is completed if we insert (\ref{theo4eqn2}) into (\ref{theo4eqn1}) and use Lemma
\ref{lemm3}.

\begin{remark} From Theorem \ref{theo4}, if~~$\|p^i-p^i_H\|_{0,\ome_s}=O(H^2),~1\leq i\leq n$, then $\|\nabla(\phi_h-\phi_h^*)\|_{0,\ome_s}=O(h)$
when $H=\sqrt{h}$, which implies Algorithm \ref{algo2} is efficient for the electrostatic potential as Algorithm \ref{algo1}. From Theorem \ref{theo5}, the two-grid solution for the density $p_h^{i,*}$ in Algorithm \ref{algo2} can achieve the same convergence rate as the finite element solution in Algorithm \ref{algo0} only when $H=h$ (cf. (\ref{phi1os})). However, the numerical results in Section 4 show the optimal convergence rate even if $H\neq h$.
\end{remark}

There is also another semi-decoupling algorithm like Algorithm \ref{algo1} as follows
%%%%%%%%%%%%%%%%%%%%%%%%
%%%%%%%%%%%%%%%%%%%%%%%% ÒÔÏ ˫ÐÐ %% ¼Ðס²¿·ÖÊÇмӽøÈ¥µÄËã·¨ 3.4 ¼°ÆäÎó²î¹À¼Æ½á¹û
\begin{algo}\label{algo3}(Two-grid algorithm III)$     $\\
Step 1. Solve the coupled problem on the coarse grid: find $p_H^i\in
S^H(\ome_s)\cap V ~(1\leq i\leq n),~~\phi_H\in S_0^H(\ome)$, such that \beas
\label{alg31}a_1(p_H^i,v_H)+b_1(p_H^i,\phi_H,v_H)=0,~~\forall
v_H\in S^H(\ome_s)\cap V_0,~~1\leq i\leq n,\eeas
\beas\label{alg32}a_2(\phi_H,w_H)+b_2(p_H^i,w_H)=\hat{f}(w_H),~~\forall v_H\in S_0^H(\ome).\eeas
Step 2. We first solve the Nernst-Planck equation on the fine grid: find
$p_h^{i,*}\in S^h(\ome_s)\cap V~(1\leq i\leq n)$, such
that \beas
\label{alg33}a_1(p_h^{i,*},v_h)+b_1(p_h^{i,*},\phi_H,v_h)=0,~~\forall
v_h\in S^h(\ome_s)\cap V_0,~~1\leq i\leq n,\eeas
then we solve the Poisson equation on the fine grid: find $\phi_h^*\in S_0^h(\ome)$, such that \beas
\label{alg34}a_2(\phi_h^*,w_h)+b_2(p_h^{i,*},w_h)=\hat{f}(w_h),~~\forall w_h\in S_0^h(\ome).\eeas
\end{algo}
In this algorithm, only the potential ${\phi _H}$ (not the density $p_H^i$ compared with Algorithm \ref{algo1}) on the coarse grid is used in Step 2. Since the error analysis and computational efficiency are similar to Algorithm \ref{algo2}, we shall not cover those results in this paper for simplicity.

%%%%%%%%%%%%%%%%%%%%%%%%%%%%%%%%%%%%%%%%%%%%%%%%%%%%%%%%%%%%%%%%%%%%%%%%%%%%%%
\section{Numerical Results}
In this section, we shall apply the two-grid algorithms to the steady-state PNP equations to verify
several derived error estimations and illustrate the efficiency of the algorithms we proposed.
\begin{exam} First we consider the steady-state PNP equations with smooth solutions as follows (see \cite{susu16} for the time-dependent PNP equations):
\beas\label{exam1} \left\{\begin{array}{ll}\nabla\cdot
\left(\nabla
p^i+q^ip^i\nabla\phi\right)=f_i,~~\mbox{in}~~ \ome,~~i=1,2,\\
-\triangle\phi-\displaystyle{\sum_{i=1}^2}
q^ip^i=f_3,~~\mbox{in}~~\ome.\end{array}\right.\eeas
Here the computational domain $\ome=[0,1]^3\subset R^3$ and $q^1=1,~~q^2=-1$.
The boundary condition and the right-hand side functions are chosen so that the exact solution $(\phi,p^1,p^2)$
is given by
\beas\left\{\begin{array}{ll}\phi=\mbox{sin}\pi x ~\mbox{sin}\pi y ~\mbox{sin}\pi z,\\
p^1=\mbox{sin}2\pi x~\mbox{sin}2\pi y ~\mbox{sin}2\pi z,\\
p^2=\mbox{sin}3\pi x~\mbox{sin}3\pi y ~\mbox{sin}3\pi z.\end{array}\right.\eeas
\end{exam}
To implement the algorithms, the code is written in Fortran and the computation is carried out on a microcomputer.
we use piecewise linear finite elements on the tetrahedral mesh to discretize the equations. For comparison of the approximation accuracy,
the finite element solution of Algorithm \ref{algo0} is solved by the following Gummel
iteration: given the initial value $(p^{1,0},p^{2,0})\in S_0^h(\ome)\times S_0^h(\ome)$, for $m\geq 0$, find $(p^{1,m+1},p^{2,m+1},\phi^{m+1})\in S_0^h(\ome)\times S_0^h(\ome) \times S_0^h(\ome)$
such that
\beas\label{giter2}(\nabla\phi^{m+1},\nabla w)+\sum_{i=1}^2q^i(p^{i,m},w)=(f_3,w),~~\forall
w\in S_0^h(\ome).\eeas
\beas\label{giter1}(\nabla p^{i,{m+1}},\nabla v)+(q^ip^{i,{m+1}}\nabla\phi^{m+1},\nabla v)=(f_i,v),~~\forall v\in
S_0^h(\ome),~~i=1,2.\eeas
The stopping criterion for this iteration is $\|\phi^{m+1}-\phi^m\|_0<10^{-5}$. The numerical results in Table 1 and 2 show that the convergence orders in $H^1$ norm and $L^2$ norm approximate first order and second order, respectively. The numerical results coincide with the theoretical results (see (\ref{phi1os}) and (\ref{errp})).

To derive the two-grid solution of Algorithm \ref{algo1}, the above Gummel iteration is used on the coarse grid and then the decoupled system is solved on the fine grid by using the coarse grid solution, see the following steps:\\
Step 1. Given the initial value $(p^{1,0},p^{2,0})\in S_0^H(\ome)\times S_0^H(\ome) $, for $m\geq 0$ find $(p^{1,m+1},p^{2,m+1},\phi^{m+1})\in S_0^H(\ome)\times S_0^H(\ome) \times S_0^H(\ome)$ such that
\bea\label{gs2}(\nabla\phi^{m+1},\nabla w)+\sum_{i=1}^2q^i(p^{i,m},w)=(f_3,w),~~\forall
w\in S_0^H(\ome).\eea
\bea\label{gs1}
(\nabla p^{i,{m+1}},\nabla v)+(q^ip^{i,{m+1}}\nabla\phi^{m+1},\nabla v)=(f_i,v),~~\forall v\in
S_0^H(\ome),~~i=1,2.\eea
Suppose $p_H^i,~i=1,2$ is the final solution of the concentration in the above iteration.\\
Step 2. We first solve the Poisson equation on the fine grid: find $\phi_h^*\in
S_0^h(\ome)$, such that
\bea\label{gs3}(\nabla \phi_h^*,\nabla w_h)+\sum_{i=1}^2q^i(p_H^{i},w_h)=(f_3,w_h),~~\forall w_h\in S_0^h(\ome),\eea
then we solve the
Nernst-Planck equation on the fine grid: find $p_h^{i,*}\in
S^h_0(\ome),~~\phi_h^*\in S_0^h(\ome)$, such that \bea\label{gs4}
(\nabla p_h^{i,*},\nabla v_h)+\sum_{i=1}^2q^i(p_h^{i,*}\nabla\phi_h^*,\nabla v_h)=(f_i,v_h),~~\forall v\in
S_0^h(\ome),~~i=1,2.\eea

The errors between the exact solution and the two-grid solutions of Algorithm \ref{algo1} are shown in Table 3 with varying mesh size and $h=H^2$. The errors in $H^1$ norm approximate the second-order reduction as $h$ becomes smaller, which indicates two-grid Algorithm \ref{algo1} remains the same order of accuracy as Algorithm \ref{algo0} but requires much less computational time than Algorithm \ref{algo0} by comparing Table 1 with Table 3.

\begin{table}[htbp]
\centering
\begin{tabular}{|*{5}{c|}}
\hline h & $\|\phi_h-\phi\|_1$ & $\|p_h^{1}-p^1\|_1$&$\|p_h^{2}-p^2\|_1$& CPU(S)\\
\hline
 1/4& 9.14E-01&3.03E+00&5.39E+00& 1.5\\ \hline
 1/8&4.80E-01&1.82E+00&3.75E+00& $-$\\ \hline
 1/16  &2.43E-01&9.57E-01&2.10E+00&7.4\\ \hline
  1/32  &1.22E-01&4.85E-01&1.09E+00&$-$\\ \hline
 1/64& 6.09E-02&2.44E-01&5.47E-01& 2933\\ \hline
\end{tabular}
\caption{The $H^1$ norm errors between the exact solution and the finite element solutions of Algorithm \ref{algo0}.
%(Example \ref{ch4ex1}) \label{ch4ex1:t1}
}
\end{table}

\begin{table}[htbp]
\centering
\begin{tabular}{|*{4}{c|}}
\hline h & $\|\phi_h-\phi\|_0$ & $\|p_h^{1}-p^1\|_0$&$\|p_h^{2}-p^2\|_0$\\
\hline
 1/4& 8.97E-02&2.41E-01&3.26E-01\\ \hline
 1/8&2.50E-02&8.99E-02&1.72E-01\\ \hline
 1/16  &6.44E-03&2.53E-02&5.59E-02\\ \hline
 1/32  &1.62E-03&6.51E-03&1.50E-02\\ \hline
 1/64& 4.06E-04&1.64E-03&3.83E-03\\ \hline
\end{tabular}
\caption{The $L^2$ norm errors between the exact solution and the finite element solutions of Algorithm \ref{algo0}.
%(Example \ref{ch4ex1}) \label{ch4ex1:t1}
}
\end{table}
\begin{table}[htbp]
\centering\label{table1}
\begin{tabular}{|*{6}{c|}}
\hline H&h=$H^2$ & $\|\phi_h^*-\phi\|_1$ & $\|p_h^{1,*}-p^1\|_1$&$\|p_h^{2,*}-p^2\|_1$& CPU(S)\\
\hline
 1/2& 1/4& 9.15E-01&3.03E+00&5.39E+00& 1.2\\ \hline
 1/4& 1/16  &2.44E-01&9.57E-01&2.10E+00&2.2\\ \hline
 1/8 & 1/64& 6.22E-02&2.44E-01&5.47E-01& 830\\ \hline
\end{tabular}
\caption{The $H^1$ norm errors between the exact solution and the two-grid finite element solutions of Algorithm \ref{algo1}.
%(Example \ref{ch4ex1}) \label{ch4ex1:t1}
}
\end{table}

The two-grid solutions of Algorithm \ref{algo2} are obtained by using the similar computational procedure in (\ref{gs2})-(\ref{gs4}), but (\ref{gs4}) is replaced with the following equation:
\beas
(\nabla p_h^{i,*},\nabla v_h)+\sum_{i=1}^2q^i(p_h^{i,*}\nabla\phi_H,\nabla v_h)=(f_i,v_h),~~\forall v\in
S_0^h(\ome),~~i=1,2.\eeas
The errors between the exact solution and the two-grid solutions of Algorithm \ref{algo2} are shown in Table 4.  The errors for the potential $\phi_h^*$ in $H^1$ norm approximate the second-order $O(H^2)$ reduction as $H$ becomes smaller,
but the error of the solutions $p_h^{1,*},~p_h^{2,*}$ can not achieve the second-order reduction, when $h=H^2$ (Comparing Table 4 with Table 1). Such an order reduction may be caused by the low order approximation for potential $\phi_H$ on coarse grid to $\phi_h$ on fine grid (See Theorem \ref{theo5}). This problem can be solved by using a smaller size of coarse grid (see the result of last row in Table 4). That means two-grid Algorithm 3.3 can also achieve the same order of accuracy as the finite element method could do, if we choose a suitable coarse mesh size $H$.

\begin{table}[htbp]
\centering\label{table1}
\begin{tabular}{|*{6}{c|}}
\hline H&h & $\|\phi_h^*-\phi\|_1$ & $\|p_h^{1,*}-p^1\|_1$&$\|p_h^{2,*}-p^2\|_1$& CPU(S)\\
\hline
 1/2& 1/4& 9.15E-01  &3.03E+00&5.39E+00& 1.2\\ \hline
 1/4& 1/16  &2.44E-01&9.89E-01&2.10E+00&2.2\\ \hline
 1/8 & 1/64& 6.22E-02&2.92E-01&5.70E-01& 830\\ \hline
 1/32 &1/64& 6.09E-02&2.46E-01&5.48E-01&1189\\ \hline
\end{tabular}
\caption{The $H^1$ norm errors between the exact solution and the two-grid finite element solutions of Algorithm \ref{algo2}.
}
\end{table}

\begin{exam}The second example is to solve (\ref{prob1})-(\ref{bd1}) for a dialanine system within a spherical computational domain with radius of $400{\AA}$. We compute the potential and concentrations of the system,
the molecular surface of which is schematically illustrated in
Fig. \ref{fig:diala}.
\end{exam}

\begin{figure}[htbp]
  \subfigure[]{
    \includegraphics[width=4.5cm,totalheight=4.5cm]{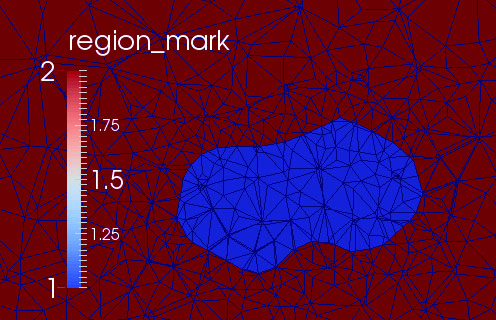}}
    \centering
     \subfigure[]{
    \includegraphics[width = 3.5cm, totalheight = 4.5cm]{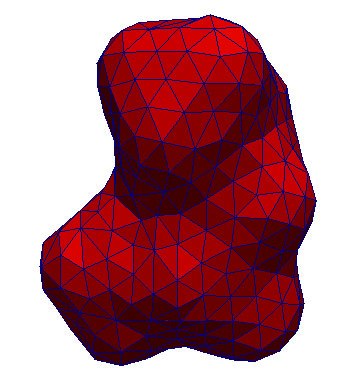}}
  \caption{An example of mesh generation for the dialanine system. (a) A close-up view of the fine mesh around the molecule. region\_mark = 1 denotes the solute region, region\_mark = 2 denotes the solvent region. (b) The triangular boundary mesh conforming to the molecular surface.}
\label{fig:diala}
\end{figure}

This system carries out a fixed total charge of $0~e_c$.
The dielectric coefficient $\epsilon$ in (\ref{prob1}) is a spacial dependent coefficient with $\epsilon_m = 2$ in
the solute region and $\epsilon_s = 78$ in the solvent region.
 We assume that there are only
monovalent ions in the salt and bulk values are set as Dirichlet
conditions in the diffusion domain. We set the
bulk densities in this 1:1 salt solution are $50 mM$. Fig.
\ref{fig:diala_var} shows the distribution of the electrostatic
potential and ion concentrations of the initial mesh.
\begin{figure}[htbp]
  \subfigure[]{
    \includegraphics[width=4cm,totalheight=4.5cm]{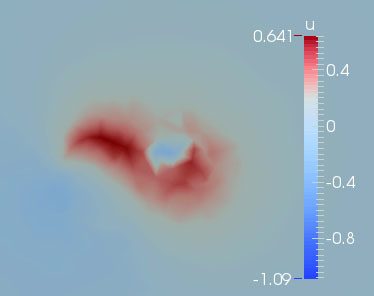}}
    \centering
     \subfigure[]{
     \includegraphics[width = 4cm, totalheight = 4.5cm]{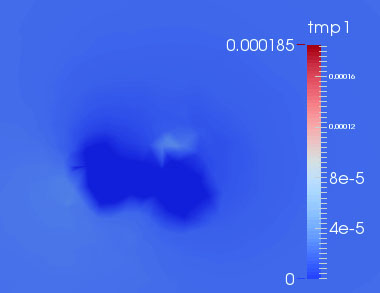}}
   \centering
    \subfigure[]{
    \includegraphics[width = 4cm, totalheight = 4.5cm]{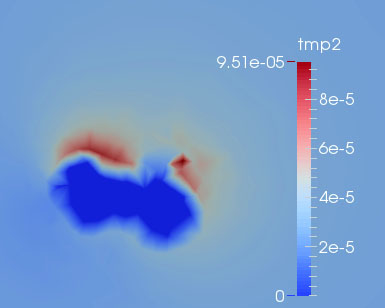}}
  \caption{Electrostatic potential and ion density(M) around the dialanine system. (a) The potential ranges from $-1.09$ to $0.64$ $kcal/mol\cdot e_c$. (b) Cation concentration distribution ranges from $0$ to $1.85\times 10^{-4}$ (M). (c) Anion concentration distribution ranges from $0$ to $9.51\times 10^{-5}$(M).}
\label{fig:diala_var}
\end{figure}

The implementation of the algorithms is based on the parallel
adaptive finite element package PHG. The parallel code is written in
C and uses MPI for message passing. The computation is carried
out on the cluster LSSC-III of the State Key Laboratory of
Scientific and Engineering Computing of China, which consists of
compute nodes with dual Intel Xeon X5550 quad-core CPUs,
interconnected via DDR InfiniBand network.

Since the surface of the molecule is extremely irregular in
practice, the initial mesh for the computation of PNP equations
should be nonuniform with a fine mesh around the surface and a
coarse mesh away from the interface. Thus in the following
experiment, we illustrate the convergence rate in terms of the order
of the degree of freedoms $N$ instead of that of mesh size
$H$. For example, if the result is $\sim\mathcal{O}(H)$, then we
think it is consistent with the estimation
$\sim\mathcal{O}(N^{-\frac{1}{3}})$ for this three dimensional problem. If we plot the log-log figure
for the original outputs (x-axis denotes the number of N, y-axis
denotes the $L^2$ or $H^1$ norm of errors), then the decay rate of the
line should be $-\frac{1}{3}$. Similarly, the decay rate should be
$-\frac{2}{3}$ if the result is $\sim\mathcal{O}(H^2)$. To estimate
the convergence rate, we refine the initial mesh step by step and
use them as coarse meshes respectively.

The numerical results for Algorithm \ref{algo1} and \ref{algo2} are shown in Fig. \ref{al1-res} and \ref{al2-res}, respectively. %the error estimations are
It is shown from Fig. \ref{al1-res} that
the convergence curve of $\displaystyle{\sum_{i}}\|p^i_h-p^{i,*}_h\|_{1,\ome_s}$ approximates
to the line with slope $-\frac 2 3$. This means the two-grid solution of the concentration in Algorithm \ref{algo1} has the optimal convergence rate which coincides with our theory in Section 3 (If $\|p^i-p_h^i\|_{0,\ome_s}=O(h^2)$, then $\displaystyle{\sum_{i}}\|p^i_h-p^i_H\|_{0,\ome_s}=O(H^2)$ and from Theorem \ref{theo3} $\displaystyle{\sum_{i}}\|p^i_h-p^{i,*}_h\|_{1,\ome_s}=O(H^2)$). Moreover, the convergence curve of $H^1$ error for the electrostatic potential decays faster than the line with slope $-\frac 2 3$ presenting a superconvergence phenomenon (see (b) in Fig. 4).
Similar results can be observed from Fig. \ref{al2-res} for Algorithm \ref{algo2}.  The concentration presents an optimal convergence phenomenon in $H^1$ norm which indicates the theoretical result in Theorem \ref{theo5} is not optimal. Similar to the result for Algorithm \ref{algo1}, the two-grid solution of electrostatic potential produces a numerically superconvergence phenomenon which may be caused by the good discrete mesh for this problem and requires further investigation.

%\begin{figure}[htbp]
%\centering
%\begin{tabular}{cc}
%\includegraphics[width=0.5\textwidth]{con_algorithm32.eps} & \includegraphics[width=0.5\textwidth]{poten_algorithm32.eps}\\
%(a)&(b)\\
%\end{tabular}
%\caption{The error estimations of Algorithm~\ref{algo1}.  (a) The convergence
%curves of $\displaystyle{\sum_i}\|p_h^i-p_h^{i,*}\|_1$; (b) The convergence curves
%of $\|\phi_h-\phi_h^*\|_1$. (DOFs represents the degree of freedoms on the coarse grid.)} \label{al1-res}
%\end{figure}
%\begin{figure}[htbp]
%\centering
%\begin{tabular}{cc}
%\includegraphics[width=0.5\textwidth]{con_algorithm33.eps} & \includegraphics[width=0.5\textwidth]{poten_algorithm33.eps}\\
%(a)&(b)\\
%\end{tabular}
%\caption{The error estimations of Algorithm~\ref{algo2}.  (a) The convergence
%curves of $\displaystyle{\sum_i}\|p_h^i-p_h^{i,*}\|_1$; (b) The convergence curves
%of $\|\phi_h-\phi_h^*\|_1$. (DOFs represents the degree of freedoms on the coarse grid.)} \label{al2-res}
%\end{figure}

% Ð޸ĴúÂë
\begin{figure}[htbp]
\centering
\begin{tabular}{cc}
\includegraphics[width=0.5\textwidth]{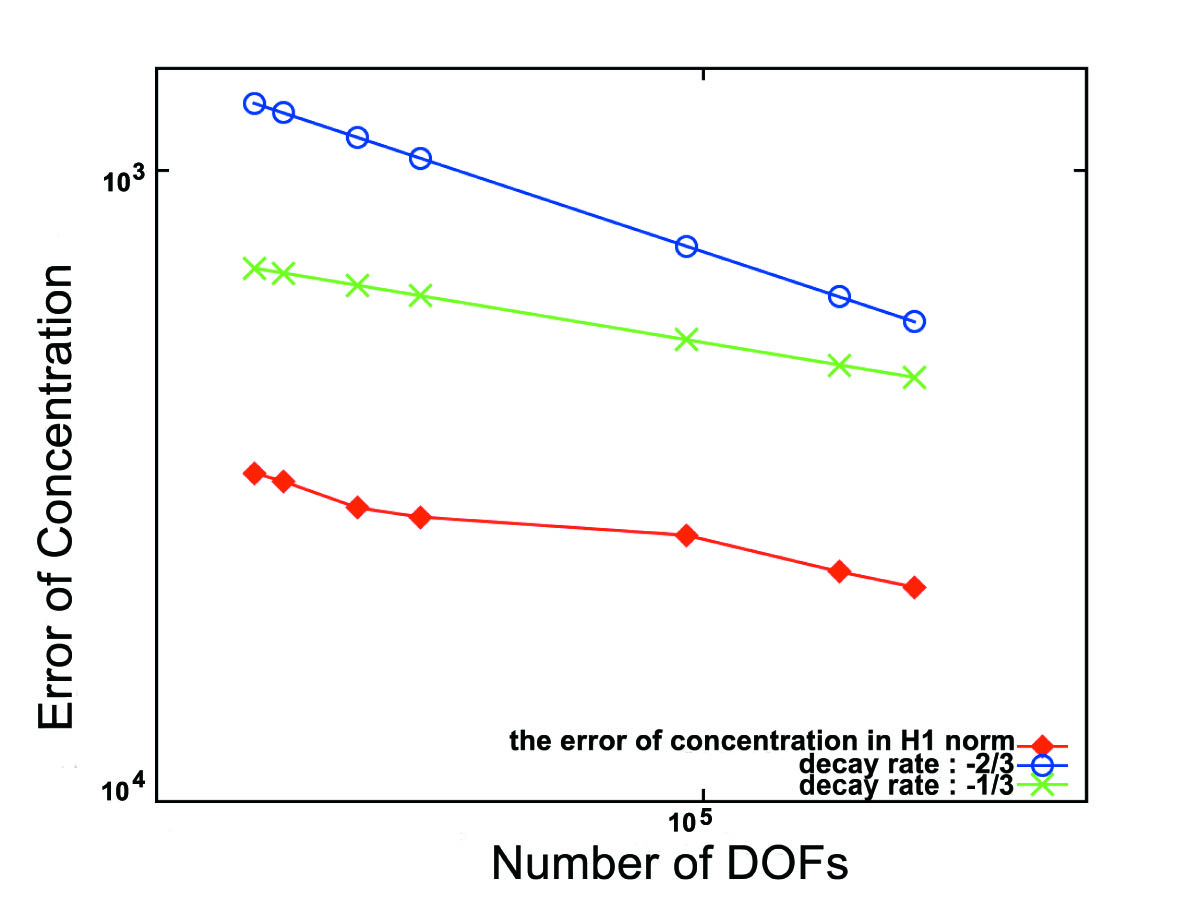} & \includegraphics[width=0.5\textwidth]{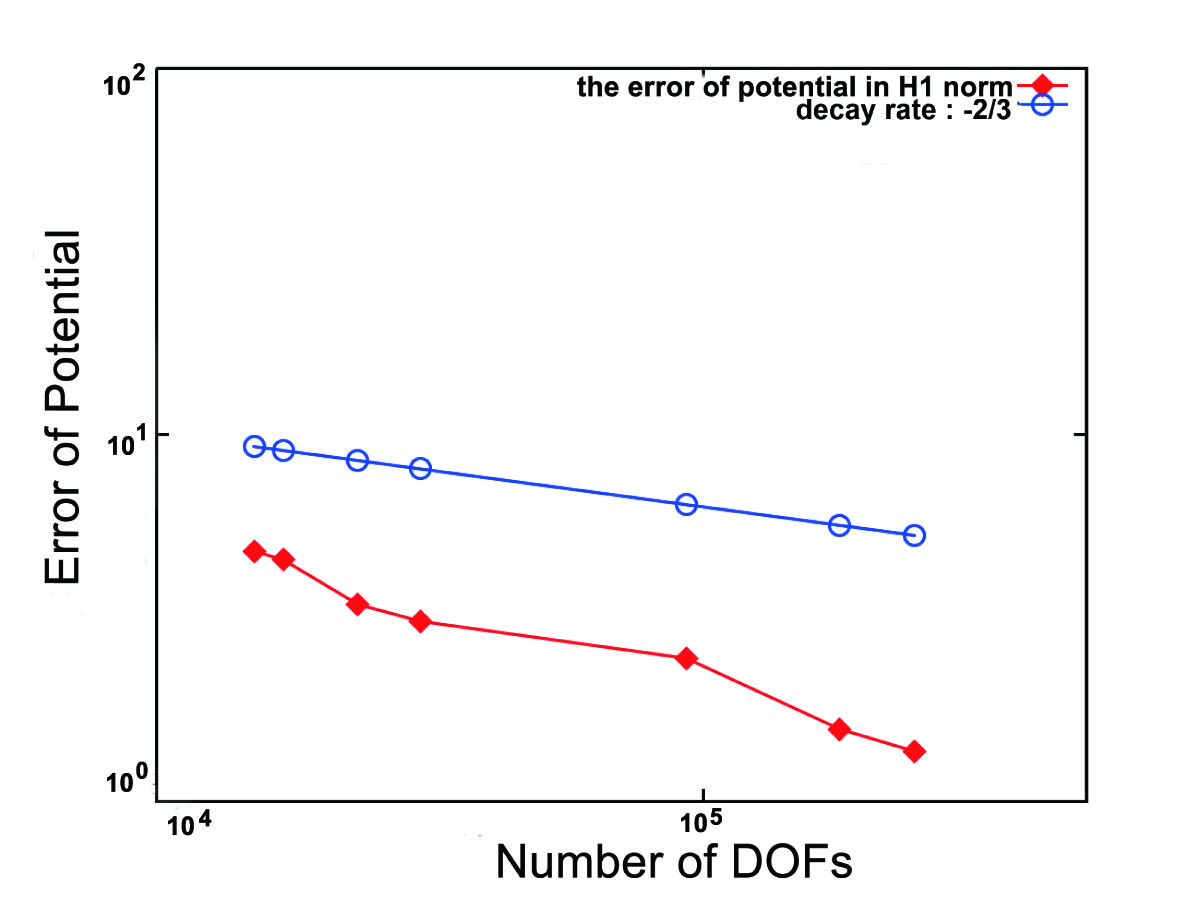}\\
(a)&(b)\\
\end{tabular}
\caption{The error estimations of Algorithm~\ref{algo1}.  (a) The convergence
curves of $\displaystyle{\sum_i}\|p_h^i-p_h^{i,*}\|_1$; (b) The convergence curves
of $\|\phi_h-\phi_h^*\|_1$. (DOFs represents the degree of freedoms on the coarse grid.)} \label{al1-res}
\end{figure}
\begin{figure}[htbp]
\centering
\begin{tabular}{cc}
\includegraphics[width=0.5\textwidth]{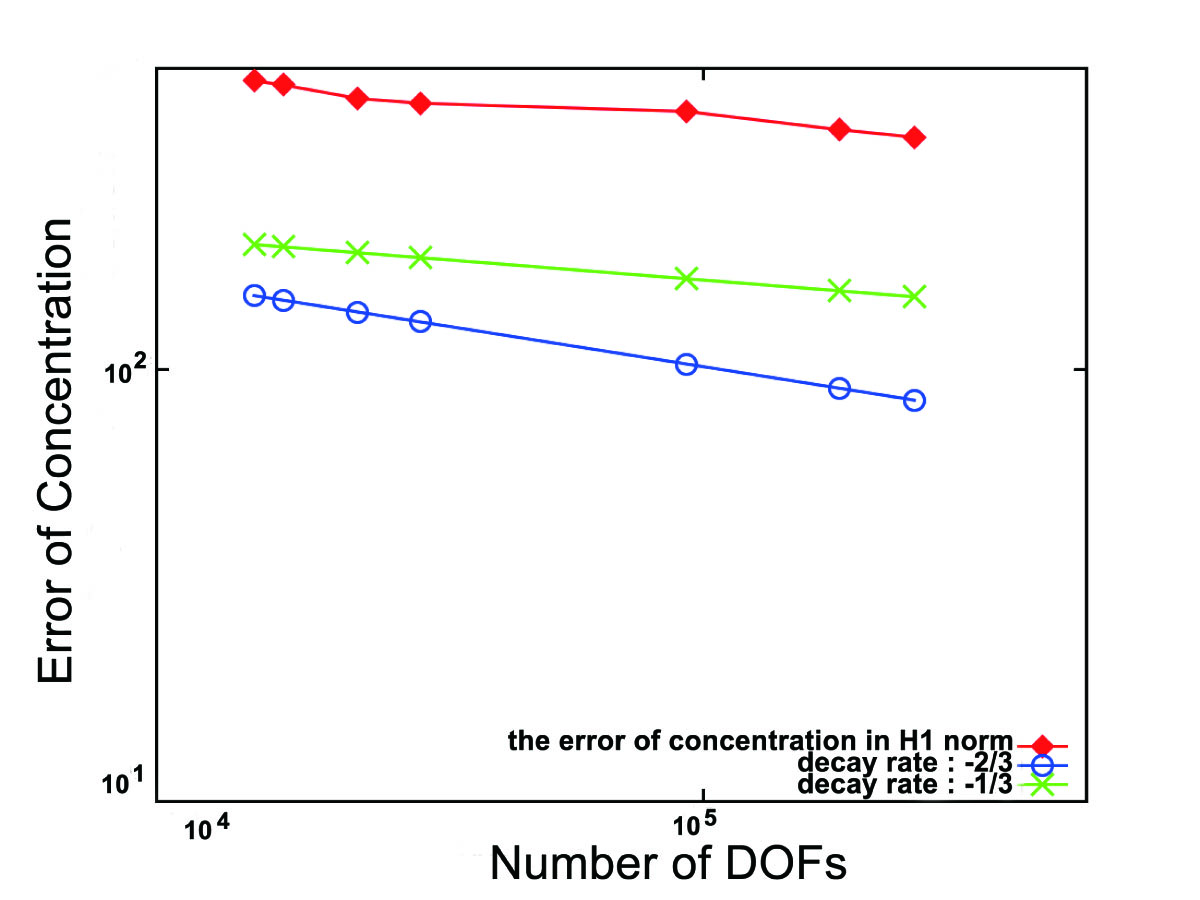} & \includegraphics[width=0.5\textwidth]{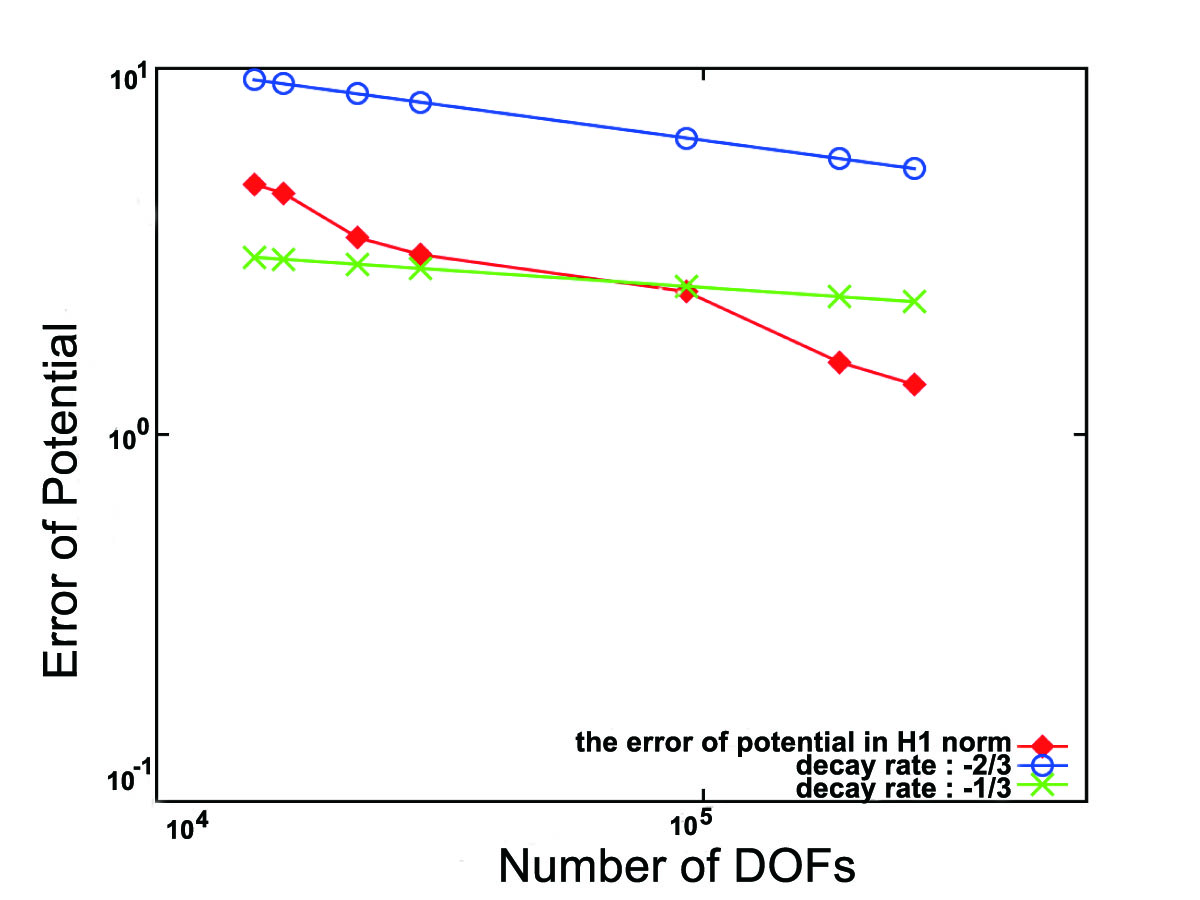}\\
(a)&(b)\\
\end{tabular}
\caption{The error estimations of Algorithm~\ref{algo2}.  (a) The convergence
curves of $\displaystyle{\sum_i}\|p_h^i-p_h^{i,*}\|_1$; (b) The convergence curves
of $\|\phi_h-\phi_h^*\|_1$. (DOFs represents the degree of freedoms on the coarse grid.)} \label{al2-res}
\end{figure}

% Ô­´úÂë
%\begin{figure}[htbp]
%\centering
%\begin{tabular}{cc}
%\includegraphics[width=0.5\textwidth]{con_algorithm3.2.eps} & \includegraphics[width=0.5\textwidth]{poten_algorithm3.2.eps}\\
%(a)&(b)\\
%\end{tabular}
%\caption{The error estimations of Algorithm~\ref{algo1}.  (a) The convergence
%curves of $\displaystyle{\sum_i}\|p_h^i-p_h^{i,*}\|_1$; (b) The convergence curves
%of $\|\phi_h-\phi_h^*\|_1$. (DOFs represents the degree of freedoms on the coarse grid.)} \label{al1-res}
%\end{figure}
%\begin{figure}[htbp]
%\centering
%\begin{tabular}{cc}
%\includegraphics[width=0.5\textwidth]{con_algorithm3.3.eps} & \includegraphics[width=0.5\textwidth]{poten_algorithm3.3.eps}\\
%(a)&(b)\\
%\end{tabular}
%\caption{The error estimations of Algorithm~\ref{algo2}.  (a) The convergence
%curves of $\displaystyle{\sum_i}\|p_h^i-p_h^{i,*}\|_1$; (b) The convergence curves
%of $\|\phi_h-\phi_h^*\|_1$. (DOFs represents the degree of freedoms on the coarse grid.)} \label{al2-res}
%\end{figure}

%%%%%%%%%%%%%%%%%%%%%%%%%%%%%%%%%%%%%%%%%%%%%%%%%%%%%%%%%%%%%%%%%%%%%%%%%%%%%%%%5
%%---------------------------------
\section{Conclusion}
In this paper, two decoupling two-grid finite element algorithms are proposed for the PNP equations. Theoretical analysis and numerical experiments show that the two-grid algorithms remain the same order of accuracy but require much less computational time compared with the classic finite element method combing with the Gummel iteration. It is promising to extend these approaches to more general settings, such as time-dependent PNP equations for ion channels, PNP equations for semiconductor devices, as well as modified PNP equations with size effects. It is also possible to generalize the framework to multilevel methods.

%%------------------------------
{\sc Acknowledgement.}
Many thanks must be expressed to Professor
Aihui Zhou and Jinchao Xu for their valuable suggestions. Thanks
also go to JingJie Xu for his discussion on the numerical experiments.
Y. Yang was supported by the China NSF (NSFC 11561016, NSFC 11561015)
and the fund from Education Department of Guangxi Province under
grant (2014GXNSFAA118004, 2014GXNSFAA118012). B. Z. Lu was supported by the
National Center for Mathematics and Interdisciplinary Sciences,
Chinese Academy of Sciences, Science Challenge Program (SCP) and the China NSF (NSFC 91530102, NSFC 21573274).
%%%%%%%%%%%%%%%%%%%%%%%%%%%%%%%%%%%%%%%%%%%%%%%%%%%%%%%%%%%%%%%%%%%%%%

\end{document}